\documentclass[a4paper,11pt,twoside,reqno]{amsart}

\usepackage[utf8]{inputenc}
\usepackage[plainpages=false,pdfpagelabels=true]{hyperref}
\usepackage{amssymb,amsthm,wasysym}
\usepackage{slashed}
\usepackage[margin=1in]{geometry}

\numberwithin{equation}{section}

\newtheorem{Satz}{Theorem}[section]
\newtheorem{Prop}[Satz]{Proposition}
\newtheorem{Lem}[Satz]{Lemma}

\theoremstyle{definition}
\newtheorem{Dfn}[Satz]{Definition}
\newtheorem{Bem}[Satz]{Remark}
\newcommand{\hess}{{\operatorname{Hess}~}}

\newcommand{\dv}{\text{dvol}_g}

\allowdisplaybreaks[1]

\renewcommand{\epsilon}{\varepsilon}

\newcommand{\R}{\ensuremath{\mathbb{R}}}
\newcommand{\C}{\ensuremath{\mathbb{C}}}
\newcommand{\D}{\slashed{D}}
\newcommand{\p}{\slashed{\partial}}

\title{Dirac-harmonic maps with potential}
\author{Volker Branding}
\date{\today}
\address{University of Vienna, Faculty of Mathematics\\
Oskar-Morgenstern-Platz 1, 1090 Vienna, Austria\\}
\email{volker.branding@univie.ac.at}
\thanks{The author gratefully acknowledges the support of the Austrian Science Fund (FWF) through the START-project Y963-N35 of Michael Eichmair and 
the project P30749-N35 “Geometric variational problems from string theory”.
}

\begin{document}

\begin{abstract}
We study the influence of an additional scalar potential on various geometric and analytic properties of Dirac-harmonic maps.
We will create a mathematical wish list of the possible benefits from inducing the potential term
and point out that the latter cannot be achieved in general.
Finally, we focus on several potentials that are motivated from supersymmetric quantum field theory.
\end{abstract} 

\maketitle

\section{Introduction and Results}
The supersymmetric nonlinear sigma model has received a lot of interest in modern
quantum field theory, in particular in string theory, over the past decades.
At the heart of this model is an action functional whose precise structure
is fixed by the invariance under various symmetry operations.

In the physics literature the model is most often formulated in the language of supergeometry
which is necessary to obtain the invariance of the action functional under supersymmetry transformations.
For the physics background of the model we refer to \cite{MR626710} and \cite[Chapter 3.4]{MR1701600}.

However, if one drops the invariance under supersymmetry transformations the resulting action functional can 
be investigated within the framework of geometric variational problems and many
substantial results in this area of research could be achieved in the last years.

To formulate this mathematical version of the supersymmetric nonlinear sigma model 
one fixes a Riemannian spin manifold \((M,g)\), a second Riemannian manifold \((N,h)\)
and considers a map \(\phi\colon M\to N\). The central ingredients in the resulting action functional
are the Dirichlet energy for the map \(\phi\) and the Dirac action
for a spinor defined along the map \(\phi\).
In the physics literature this spinor would also take values in a Grassmann algebra turning it
into a non-commuting object.
Although we neglect the invariance under supersymmetry transformations 
the action functional is invariant under diffeomorphisms on the domain
and for a two-dimensional domain it is also invariant under conformal transformations of the metric on
the domain. The first invariance gives rise to the stress-energy tensor which is 
conserved whenever we consider a critical point of the action functional.
The conformal invariance in dimension two leads to various nice properties
of the critical points such as the possibility of removing isolated point singularities
whenever a certain energy is finite.

The mathematical study of the supersymmetric nonlinear sigma model 
with standard spinors was initiated in \cite{MR2262709}
where the terminology \emph{Dirac-harmonic maps} for the critical points
was established. The Dirac-harmonic map equations consist of a
semilinear second order elliptic equation for the map \(\phi\)
and a linear Dirac equation, which is elliptic and of first order, for the spinor along the map.

Motivated from various variants in the physics literature several
extensions of Dirac-harmonic maps have also been studied from a mathematical point of view:
Taking into account a two-form contribution in the action functional
one is led to \emph{magnetic Dirac-harmonic maps} \cite{MR3305429},
Dirac-harmonic maps to target spaces with torsion are investigated in \cite{MR3493217}.
Adding a curvature term to the action functional, which is quartic in the spinors,
one obtains \emph{Dirac-harmonic maps with curvature term},
which have been studied extensively in \cite{MR3333092,MR3558358,MR3735550,MR3830780,MR4034775,MR2370260,MR4018319}.
Recently, another extension of Dirac-harmonic maps receives growing interest:
Here, one considers an additional field in the action functional, the so-called \emph{gravitino} \cite{MR3772035}.
In the physics literature the gravitino is the supersymmetric partner of the metric on the domain
meaning that both can be transformed into each other by a supersymmetry transformation.

At present, many results on the qualitative behavior of a given Dirac-harmonic map
are known as for example the regularity of weak solutions \cite{MR2544729}.

However, it remains a challenging mathematical problem to prove a general existence result for Dirac-harmonic maps.
The existence of uncoupled Dirac-harmonic maps using the Atiyah-Singer index theorem was established in \cite{MR3070562}.
The terminology uncoupled here refers to the fact that these Dirac-harmonic maps are constructed 
out of a given harmonic map such that all terms in the Euler-Lagrange equations for Dirac-harmonic maps
vanish individually. An existence result for Dirac-harmonic maps from surfaces with boundary
via the heat flow method could be achieved in \cite{MR3724759}.

Let us mention another recent approach to the existence problem for Dirac-harmonic maps.
Motivated from the classical work of Sacks and Uhlenbeck \cite{MR604040} for harmonic maps,
the so-called $\alpha$\emph{-Dirac-harmonic maps} haven been introduced and investigated
in a series of articles
\cite{MR4259181,MR4310166,MR4232501}
and some first existence results for Dirac-harmonic maps from closed surfaces could be achieved within this framework.

There are also some existence results available for the hyperbolic version of Dirac-harmonic maps, which are \emph{Dirac-wave maps},
and their extensions. These arise if one considers a domain manifold with a Lorentzian metric.
Existence results for Dirac-wave maps from two-dimensional Minkowski space 
have been obtained in \cite{MR3917346,MR2138082}.
The Cauchy problem for Dirac-wave maps with curvature term on expanding spacetimes
has been successfully studied in \cite{MR3830277}.

Besides the aforementioned existence results several Liouville-type results 
have also been established \cite{MR3830780,MR3886921,MR4034775,MR2370260}. These provide criteria under which a Dirac-harmonic map
must be trivial, that is the map part maps to a point and the spinor vanishes identically.

For the current state of research on the mathematical aspects of the supersymmetric nonlinear
sigma model we refer to the recent survey article \cite{MR3913850}.

In the physics literature there exists a version of
the supersymmetric nonlinear sigma model coupled to a scalar potential.
This potential has to be chosen in such a way that it respects the invariance
of the action functional under supersymmetry transformations.
This particular potential was introduced in \cite{MR719813}, see also
the discussion in \cite[Theorem 3.82]{MR1701600}.

Besides the mathematical results that build on the nonlinear supersymmetric sigma model
let us also recall several results on harmonic maps coupled to a scalar potential
that were studied in the mathematics literature.
Harmonic maps coupled to a scalar potential were introduced in \cite{MR1433176},
the corresponding heat flow was studied in \cite{MR1800592}.
For the current status of research on harmonic maps with potential
we refer to the introduction of \cite{MR3673634},
the regularity of weak harmonic maps with (smooth) potential from surfaces
was recently established in \cite[Theorem 1.1]{MR3884770}.

In this article we will investigate the action functional for Dirac-harmonic maps
coupled to an arbitrary scalar potential.
Initially, one could hope that it is favorable to also include such a potential term
in the action functional. Many of the difficulties in the analysis of Dirac-harmonic maps
have their origin in the fact that the corresponding action functional is unbounded from below.
An additional potential term can in principle be used to repair this flaw.

Let us now describe the mathematical setup that we employ in more detail. 
We assume that \((M,g)\) is a closed Riemannian spin manifold with spinor bundle \(\Sigma M\), 
for more details about spin geometry see the book \cite{MR1031992}. Moreover, let \((N,h)\) be a second closed Riemannian
manifold. 
Let \(\phi\colon M\to N\) be a map,
integrating the square of its differential \(d\phi\in\Gamma(T^*M\otimes\phi^\ast TN)\)
leads to the usual Dirichlet energy.
Together with the pullback bundle \(\phi^\ast TN\) 
we consider the twisted bundle \(\Sigma M\otimes\phi^\ast TN\). 
The induced connection on this bundle will be denoted by \(\tilde{\nabla}\). 
Moreover, we have an induced Hermitian scalar product on \(\Sigma M\otimes\phi^\ast TN\)
of which we will always take the real part.
Sections \(\psi\in\Gamma(\Sigma M\otimes\phi^\ast TN)\) in this bundle are called \emph{vector spinors}.
The natural operator acting on vector spinors
is the twisted Dirac operator, denoted by \(\D\).  
More precisely, the twisted Dirac operator is given by \(\D=e_i\cdot\tilde{\nabla}_{e_i}\), where \(\{e_i\},i=1,\ldots,m=\dim M\) is
an orthonormal basis of \(TM\) and \(\cdot\) denotes Clifford multiplication. 
It is an elliptic, first order operator, which is self-adjoint with respect to the \(L^2\)-norm.
We are using the Einstein summation convention, that is we sum over repeated indices.
Whenever choosing indices we will use Latin letters for indices related to \(M\) and Greek letters for indices on \(N\).

Recall that Clifford multiplication is skew-symmetric, namely
\[
\langle\chi,X\cdot\xi\rangle_{\Sigma M}=-\langle X\cdot\chi,\xi\rangle_{\Sigma M}
\]
for all \(\chi,\xi\in\Gamma(\Sigma M)\) and all \(X\in TM\).
In addition, the Clifford relations
\[
X\cdot Y+Y\cdot X=-2g(X,Y)
\]
hold for all \(X,Y\in TM\).

In terms of local coordinates \(y^\alpha\) on \(N\) the vector spinor \(\psi\) can be expanded as \(\psi=\psi^\alpha\otimes\frac{\partial}{\partial y^\alpha}\), 
and thus the twisted Dirac operator \(\D\) is locally given by
\[
\D\psi=\big(\p\psi^\alpha+\Gamma^\alpha_{\beta\gamma}(\phi)\nabla\phi^\beta\cdot\psi^\gamma\big)\otimes\frac{\partial}{\partial y^\alpha}.
\]
Here, \(\p:=e_i\cdot\nabla^{\Sigma M}_{e_i}\colon\Gamma(\Sigma M)\to\Gamma(\Sigma M)\) denotes the standard Dirac operator and \(\Gamma^\alpha_{\beta\gamma}\) are the Christoffel symbols
of the manifold \(N\).

In order to define Dirac-harmonic maps with potential assume that \(V\colon \Sigma M\otimes\phi^\ast TN\times N\to\R\) is a function which we
assume to be smooth in the following.

The central object in this article is the action functional
\begin{equation}
\label{action-functional}
S_P(\phi,\psi)=\int_M\big(|d\phi|^2+\langle\psi,\D\psi\rangle-2V(\phi,\psi)\big)\dv,
\end{equation}
where we use the subscript \(P\) to highlight the presence of the potential \(V(\phi,\psi)\).

Before we turn to the mathematical analysis of the action functional \eqref{action-functional}
let us give a mathematical wish list of what the potential \(V(\phi,\psi)\) could achieve:
\begin{enumerate}
 \item The potential \(V(\phi,\psi)\) removes the indefinite character of the action functional \eqref{action-functional}
 such that \(S_P(\phi,\psi)\) is now bounded from below. Moreover, the potential can help in establishing an existence result
 for critical points of \eqref{action-functional}.
 \item The potential also ensures that critical points of \eqref{action-functional} are  stable in the sense that the
  second variation is positive under certain geometric or analytic assumptions.
 \item The regularity theory for weak solutions of the Dirac-harmonic map system (which are the critical points of \eqref{action-functional} in the
  case of a vanishing potential) also applies to the critical points of \eqref{action-functional}.
 \item The potential respects the invariance under diffeomorphisms on the domain and also under 
 conformal transformations in the case of a two-dimensional domain manifold.
\end{enumerate}

Even without going into the mathematical details one should expect that there cannot be a 
potential \(V(\phi,\psi)\) that satisfies all items of the above wish list.
For example, if we have a potential that is able to suppress the negative eigenvalues
of the twisted Dirac operator \(\D\) then we cannot expect it to be integrable in the Sobolev
spaces that are used in the regularity analysis of weak Dirac-harmonic maps.

The inclusion of a fine-tuned potential term turned out to be 
very helpful in order to derive existence results for perturbed versions
of Dirac-harmonic maps as was done in \cite{MR2860404,MR3978518,isobe2017morse}.
However, one still needs highly sophisticated techniques from the calculus of variation
in order to achieve such results.

This article is organized as follows:
In Section 2 we calculate the first and the second variation of \eqref{action-functional}.
The third section studies the possible benefits of the scalar potential \(V(\phi,\psi)\).
In particular, we investigate if it
can help to ensure the positivity of the action functional \eqref{action-functional} and its second variation.
Moreover, we make some comments on the importance of the second variation
of an action functional for quantum field theory.
Section 4 focuses on the regularity of weak Dirac-harmonic maps with potential.

The last section discusses various explicit potentials, both from the physics literature
and inspired from the mathematical wish list above,
and their possible benefit for the mathematical analysis.

\section{Variational formulas}
In this section we derive several variational formulas for Dirac-harmonic maps with potential.
In the following \(\epsilon\) will always denote a small positive real number.

\subsection{First variation}
Let us briefly discuss the first variation of the action functional \eqref{action-functional}.
\begin{Prop}
The critical points of the action functional \eqref{action-functional} are given by
\begin{align}
\label{euler-lagrange-phi}\tau(\phi)=&\frac{1}{2}R^N(\psi,e_i\cdot\psi)d\phi(e_i)-V_\phi(\phi,\psi), \\
\label{euler-lagrange-psi}\D\psi=&V_{\psi}(\phi,\psi),
\end{align}
where \(\tau(\phi)\) is the tension field of the map \(\phi\), 
\(V_{\phi}(\phi,\psi)\) is the functional derivative 
of the potential \(V(\phi,\psi)\) with respect to the map \(\phi\)
and \(V_{\psi}(\phi,\psi)\) represents the functional derivative 
of the potential \(V(\phi,\psi)\) with respect to the spinor \(\psi\).
\end{Prop}

\begin{proof}
We consider a variation of the pair \((\phi,\psi)\), given by \(\phi_t\colon (-\epsilon,\epsilon)\times M\to N\) and 
\(\psi_t\colon(-\epsilon,\epsilon)\times M\to\Sigma M\otimes\phi_t^\ast TN\)
that satisfies \(\frac{\partial\phi_t}{\partial t}\big|_{t=0}=\eta\) and \(\frac{\tilde\nabla\psi_t}{\partial t}\big|_{t=0}=\xi\).

The following variational formulas are well known, see for example \cite[Section 4]{MR3772035}:
\begin{align*}
\frac{d}{dt}\big|_{t=0}\int_M|d\phi_t|^2\dv&=-2\int_M\langle\tau(\phi),\eta\rangle\dv,\\
\frac{d}{dt}\big|_{t=0}\int_M\langle\psi_t,\D\psi_t\rangle\dv&=\int_M(2\langle\xi,\D\psi\rangle+\langle\eta,R^N(\psi,e_i\cdot\psi)d\phi(e_i)\rangle)\dv.
\end{align*}

In order to obtain the variation of the potential \(V(\phi,\psi)\) we employ the chain rule and find
\begin{align*}
\frac{d}{dt}\big|_{t=0} V(\phi_t,\psi_t)&=\langle V_{\phi_t}(\phi_t,\psi_t),d\phi_t(\partial_t)\rangle|_{t=0}
+\langle V_{\psi_t}(\phi_t,\psi_t),\frac{\tilde\nabla\psi_t}{\partial t}\rangle|_{t=0}\\
&=\langle V_\phi(\phi,\psi),\eta\rangle+\langle V_{\psi}(\phi,\psi),\xi\rangle.
\end{align*}
The claim follows from combining the different equations.
\end{proof}

Let us give the local version of the Euler-Lagrange equations \eqref{euler-lagrange-phi}, \eqref{euler-lagrange-psi}.
To this end let \((U,x^i)\) be a local chart on \(M\) and \((V,y^\alpha)\) be a local
chart on \(N\) such that \(\phi(U)\subset V\). 
Then the system \eqref{euler-lagrange-phi}, \eqref{euler-lagrange-psi} acquires the form
\begin{align*}
\Delta\phi^\alpha=&-\Gamma^\alpha_{\beta\gamma}g^{ij}\frac{\partial\phi^\beta}{\partial x^i}\frac{\partial\phi^\gamma}{\partial x^j}
+\frac{1}{2}R^\alpha_{~\beta\gamma\delta}\frac{\partial\phi^\beta}{\partial x^i}\langle\psi^\gamma,e_i\cdot\psi^\delta\rangle_{\Sigma M}
-(V_\phi(\phi,\psi))^\alpha,\\
\p\psi^\alpha=&-\Gamma^\alpha_{\beta\gamma}e_i\cdot\frac{\partial\phi^\beta}{\partial x_i}\psi^\gamma
+(V_\psi(\psi,\phi))^\alpha.
\end{align*}

\begin{Bem}
The functional derivative of \(V(\phi,\psi)\) with respect to \(\psi\) in \eqref{euler-lagrange-psi}
can be defined as follows: Consider a variation of the vector spinor \(\psi\) that is 
\(\psi_t\colon(-\epsilon,\epsilon)\times M\to\Sigma M\otimes\phi^\ast TN\)
that satisfies \(\frac{\tilde\nabla\psi_t}{\partial t}\big|_{t=0}=\xi\) while keeping the map \(\phi\) fixed.
Then \(V_\psi(\phi,\psi)\in\Gamma(\Sigma M\otimes\phi^\ast TN)\) is defined as follows
\begin{align}
\langle V_\psi(\phi,\psi),\xi\rangle:=\frac{d}{dt}\big|_{t=0}V(\phi,\psi_t).
\end{align}
\end{Bem}

\subsection{Second variation}
In this section we will calculate the second variation of the action functional \eqref{action-functional}.
In general, one should of course expect the second variation of \eqref{action-functional} to be indefinite due to the
presence of the Dirac term.
Again, we make us of a variation of the pair \((\phi,\psi)\), given by \(\phi_t\colon (-\epsilon,\epsilon)\times M\to N\) and 
\(\psi_t\colon(-\epsilon,\epsilon)\times M\to\Sigma M\otimes\phi_t^\ast TN\) satisfying
\begin{align*}
\frac{\partial\phi_t}{\partial t}\big|_{t=0}=\eta,\qquad 
\frac{\tilde\nabla\psi_t}{\partial t}\big|_{t=0}=\xi.
\end{align*}

Using the first variation \eqref{euler-lagrange-phi}, \eqref{euler-lagrange-psi} we obtain
\begin{align*}
\frac{d}{dt}S_P(\phi_t,\psi_t)=&2\int_M\big(\langle d\phi_t(\partial_t),-\tau(\phi_t)+\frac{1}{2}R^N(\psi_t,e_i\cdot\psi_t)d\phi_t(e_i)
-V_\phi(\phi_t,\psi_t)\rangle \\
&+\langle\frac{\tilde\nabla\psi_t}{\partial t},\D\psi_t-V_{\psi_t}(\phi_t,\psi_t)\rangle
\big) \dv.
\end{align*}
In order to obtain the formula for the second variation of \eqref{action-functional} we have
to differentiate the above equation with respect to \(t\) once more and evaluate it at \(t=0\).

It is well-known that 
\begin{align*}
\frac{d}{dt}\big|_{t=0}\int_M\langle d\phi_t(\partial_t),\tau(\phi_t)\rangle\dv=&\int_M\langle\frac{\partial^2\phi_t}{\partial t^2},\tau(\phi_t)\rangle\dv\big|_{t=0} \\
&+\int_M(-|\nabla\eta|^2+\langle R^N(\eta,d\phi)d\phi,\eta\rangle)\dv.
\end{align*}

Moreover, we find that
\begin{align*}
\frac{d}{dt}\big|_{t=0}\big(\frac{1}{2}R^N(\psi_t,e_i\cdot\psi_t)d\phi_t(e_i)\big)=&\frac{1}{2}(\nabla_\eta R^N)(\psi,e_i\cdot\psi)d\phi(e_i)
+R^N(\xi,e_i\cdot\psi)d\phi(e_i) \\
&+\frac{1}{2}R^N(\psi,e_i\cdot\psi)\nabla_{e_i}\eta
\end{align*}
leading to
\begin{align*}
\frac{d}{dt}\big|_{t=0}\int_M&\langle d\phi_t(\partial_t),\frac{1}{2}R^N(\psi_t,e_i\cdot\psi_t)d\phi_t(e_i)\rangle\dv\\
=&\int_M\big(\langle\frac{\partial^2\phi_t}{\partial t^2},\frac{1}{2}R^N(\psi_t,e_i\cdot\psi_t)d\phi_t(e_i)\rangle\dv\big|_{t=0}
+\frac{1}{2}\langle\eta,(\nabla_\eta R^N)(\psi,e_i\cdot\psi)d\phi(e_i)\rangle \\
&+\langle\eta,R^N(\xi,e_i\cdot\psi)d\phi(e_i)\rangle
+\frac{1}{2}\langle\eta,R^N(\psi,e_i\cdot\psi)\nabla_{e_i}\eta\rangle \big)\dv.
\end{align*}

The second variation of the Dirac action can be calculated as
\begin{align*}
\frac{d}{dt}\big|_{t=0}\int_M\langle\frac{\tilde\nabla\psi_t}{\partial t},\D\psi_t\rangle\dv=&\int_M\langle\frac{\tilde\nabla^2\psi_t}{\partial t^2},\D\psi_t\rangle\dv\big|_{t=0} \\
&+\int_M\big(\langle\eta,R^N(\xi,e_i\cdot\psi)d\phi(e_i)\rangle+\langle\xi,\D\xi\rangle\big)\dv.
\end{align*}

Regarding the potential term \(V(\phi,\psi)\) we find by using the chain rule
\begin{align*}
\frac{d^2}{dt^2}\big|_{t=0}V(\phi_t,\psi_t)
=&\frac{d}{dt}\big|_{t=0}\big(\langle V_{\phi_t}(\phi_t,\psi_t),d\phi_t(\partial_t)\rangle
+\langle V_{\psi_t}(\phi_t,\psi_t),\frac{\tilde\nabla\psi_t}{\partial t}\rangle\big)\\
=&\langle\frac{\partial^2\phi_t}{\partial t^2},V_{\phi_t}(\phi_t,\psi_t)\rangle\big|_{t=0}
+\langle\frac{\tilde{\nabla}^2\psi_t}{\partial t^2},V_\psi(\phi_t,\psi_t)\rangle\big|_{t=0} \\
&+\iota(\eta,\eta)\big(V_{\phi\phi}(\phi,\psi)\big)+\iota(\xi,\xi)\big(V_{\psi\psi}(\phi,\psi)\big)
+2\iota(\eta,\xi)\big(V_{\phi\psi}(\phi,\psi)\big).
\end{align*}

Here, and in the following, the symbol \(\iota\) is used to either denote the insertion of a vector field \(\eta\)
or a vector spinor \(\xi\) into the corresponding functional derivative of the potential.

Combining the previous steps we get
\begin{Satz}
Let \((\phi,\psi)\) be a smooth solution of the system \eqref{euler-lagrange-phi}, \eqref{euler-lagrange-psi}.
Then the second variation of the action functional \eqref{action-functional} is given by
\begin{align}
\label{second-variation}
\frac{d^2}{dt^2}\big|_{t=0}S_p(\phi_t,\psi_t)=2\int_M\bigg(
&|\nabla\eta|^2-\langle R^N(\eta,d\phi)d\phi,\eta\rangle+\langle\xi,\D\xi\rangle \\
\nonumber&+\frac{1}{2}\langle\eta,(\nabla_\eta R^N)(\psi,e_i\cdot\psi)d\phi(e_i)\rangle 
+2\langle\eta,R^N(\xi,e_i\cdot\psi)d\phi(e_i)\rangle \\
\nonumber&+\frac{1}{2}\langle\eta, R^N(\psi,e_i\cdot\psi)\nabla_{e_i}\eta\rangle
-\iota(\eta,\eta)\big(V_{\phi\phi}(\phi,\psi)\big) \\
\nonumber&-\iota(\xi,\xi)\big(V_{\psi\psi}(\phi,\psi)\big)
-2\iota(\eta,\xi)\big(V_{\phi\psi}(\phi,\psi)\big)
\bigg)\dv.
\end{align}
\end{Satz}
\begin{proof}
This follows from the previous equations and the fact that 
\((\phi,\psi)\) is a smooth solution of the system \eqref{euler-lagrange-phi}, \eqref{euler-lagrange-psi}.
\end{proof}

Formally, we can interpret the second variation of \eqref{action-functional} in the following way
\begin{align*}
\hess\big(S_P(\phi,\psi)\big)(\eta,\xi)=&\int_M\langle I(\eta,\xi),(\eta,\xi)\rangle\dv,
\end{align*}
where \(I(\phi,\psi)\colon\Gamma(\phi^\ast TN)\times\Gamma(\Sigma M\otimes\phi^\ast TN)\to\Gamma(\phi^\ast TN)\times\Gamma(\Sigma M\otimes\phi^\ast TN)\)
is a semilinear elliptic differential operator of order 2 that generalizes the well-studied \emph{Jacobi operator}.

We can define the associated Jacobi operator of \eqref{action-functional} as follows
\begin{align}
\label{jacobi-operator}
I(\eta,\xi)=2\big(&-\Delta\eta+R^N(d\phi,\eta)d\phi+\frac{1}{2}(\nabla_\eta R^N)(\psi,e_i\cdot\psi)d\phi(e_i)
+2R^N(\xi,e_i\cdot\psi)d\phi(e_i)\\
\nonumber&+\frac{1}{2}R^N(\psi,e_i\cdot\psi)\nabla_{e_i}\eta 
-\iota(\eta)(V_{\phi\phi}(\phi,\psi))-2\iota(\xi)(V_{\phi\psi}(\phi,\psi)),\\
\nonumber&\D\xi-\iota(\xi)V_{\psi,\psi}(\phi,\psi)\big).
\end{align}
Here, we have \(\eta\in\Gamma(\phi^\ast TN), \xi\in\Gamma(\Sigma M\otimes\phi^\ast TN)\)
and \(\Delta\) is the rough Laplacian on the vector bundle \(\phi^\ast TN\).
Note that the Jacobi operator cannot be defined in a unique fashion as we 
can assign \(\eta\) and \(\xi\) contained in the term
\(\langle\eta,R^N(\xi,e_i\cdot\psi)d\phi(e_i)\rangle\)
to either the equation for \(\eta\) or the equation for \(\xi\).

\section{Positivity aspects and the question of stability}
The initial mathematical motivation to study Dirac-harmonic maps coupled to a scalar
potential was the hope that a suitable potential may cure at least some of the analytic
difficulties that arise in the analysis of Dirac-harmonic maps, in particular the problems
arising from the unboundedness of the action functional. In this section we will analyze this idea
in some more detail.

First, we will have a look at the action functional \eqref{action-functional}.
If we assume that \(M\) is compact, then we know from general spectral theory 
that the twisted Dirac operator \(\D\) has a discrete spectrum
consisting of an unbounded discrete sequence of positive and negative eigenvalues.
If we expand a vector spinor \(\psi\) with respect to a basis of eigenspinors \(\psi_I\)
as
\begin{align*}
\psi:=\sum_{I=-\infty}^\infty\alpha_I\psi_I,\qquad \alpha_I\in\C,
\end{align*}
the Dirac term in the action functional \eqref{action-functional} acquires the form
\begin{align*}
\int_M\langle\psi,\D\psi\rangle\dv=\sum_{J=-\infty}^\infty|\alpha_J|^2\lambda_J.
\end{align*}
Here, we used that \(\D\psi_J=\lambda_J\psi_J\), where \(\lambda_J\) represents
the J-th eigenvalue of the twisted Dirac operator \(\D\).

Consequently, if we want to obtain a lower bound for the action functional \eqref{action-functional}
we need to require that
\begin{align}
\label{positivity-action}
\sum_{J=-\infty}^\infty|\alpha_J|^2\lambda_J-2\int_MV(\phi,\psi)\dv\geq 0.
\end{align}
We can conclude that the potential \(V(\phi,\psi)\) needs to be able to approach infinity
in the same fashion as the spectrum of the twisted Dirac operator \(\D\) does.
However, if we want to allow the potential \(V(\phi,\psi)\) to be smooth
and to choose a compact target manifold \(N\) then it will not be able
to compensate the lowest eigenvalue of the Dirac operator.
On the other hand, if we choose a non-compact target manifold \(N\)
we might encounter additional difficulties in dealing with the analytic aspects of 
the critical points of \eqref{action-functional}.

\subsection{Positivity of the second variation}
In the geometric calculus of variation the second variation of a given action functional
encodes important information on the stability of a critical point.
In the case of harmonic maps it is well known that the second variation of the Dirichlet energy
is positive when the target manifold has negative curvature
such that critical points are stable then \cite{MR214004}.

In the case of Dirac-harmonic maps (with or without additional potential)
one again cannot expect the second variation of the action functional to be positive in general
due to the presence of the Dirac term.
Let us again derive a condition on the potential \(V(\phi,\psi)\) that would
be sufficient to ensure the positivity of the second variation.
To this end we estimate the second variation \eqref{second-variation} as follows:

\begin{align*}
\hess\big(S_P(\phi,\psi)\big)(\eta,\xi)\geq2\int_M&\bigg(
|\nabla\eta|^2+\langle\xi,\D\xi\rangle \\
&-c_1|\eta|^2|d\phi|^2-c_2|\eta|^2|\psi|^2|d\phi|
-c_3|\eta||\xi||\psi||d\phi|-c_4|\eta||\psi|^2|\nabla\eta| \\
&-\iota(\eta,\eta)\big(V_{\phi\phi}(\phi,\psi)\big)
-\iota(\xi,\xi)\big(V_{\psi\psi}(\phi,\psi)\big) \\
&-2\iota(\eta,\xi)\big(V_{\phi\psi}(\phi,\psi)\big)
\bigg)\dv\\
\geq 2\int_M&\bigg(
\frac{1}{2}|\nabla\eta|^2+\langle\xi,\D\xi\rangle \\
&-c_5|\eta|^2(|d\phi|^2+|\psi|^4)-c_6|\xi|^2|\psi|^2\\
&-\iota(\eta,\eta)\big(V_{\phi\phi}(\phi,\psi)\big)
-\iota(\xi,\xi)\big(V_{\psi\psi}(\phi,\psi)\big) \\
&-2\iota(\eta,\xi)\big(V_{\phi\psi}(\phi,\psi)\big)
\bigg)\dv.
\end{align*}
Here \(c_i,i=1,\ldots,6\) are positive constants
that only depend on the geometry of \(N\).

In the analysis of Dirac-harmonic maps and their variants the following energy often
turns out to be crucial
\begin{align*}
e(\phi,\psi)=|d\phi|^2+|\psi|^4.
\end{align*}

The second variation of \eqref{action-functional} would be positive if the following conditions hold:
\begin{enumerate}
 \item The Hessian of the potential with respect to the map \(\phi\) needs to satisfy
 \begin{align}
  \label{second-variation-positivity-a}
  -\iota(\eta,\eta)\big(V_{\phi\phi}(\phi,\psi)\big)\geq C|\eta|^2e(\phi,\psi)
  \end{align}
  for some positive constant \(C\).
  Again, we realize that we would need to choose a non-compact target manifold
  in order for the Hessian of the potential to be negative globally.

 \item In addition, we need to require the positivity of the following expression
 \begin{align}
 \label{second-variation-positivity-b}
\int_M\bigg(\langle\xi,\D\xi\rangle-c_6|\xi|^2|\psi|^2
-\iota(\xi,\xi)\big(V_{\psi\psi}(\phi,\psi)\big) 
-2\iota(\eta,\xi)\big(V_{\phi\psi}(\phi,\psi)\big)\bigg)\dv\geq 0.
 \end{align}
\end{enumerate}
The analysis carried out in this section strongly suggests that it will not be
possible to find a ``suitable'' potential \(V(\phi,\psi)\) that
can ensure the positivity of the second variation of the action functional \eqref{action-functional}
and that also has a reasonable application in quantum field theory.

\subsection{Some remarks on the notion of stability and quantum field theory}
Let us also make a remark on why the stability of a critical point of an action functional,
within the framework of the geometric calculus of variations,
could have important implications for calculations in quantum field theory.

In the context of geometric variational problems the second variation of an action functional 
gives important information on the stability of a critical point.
Roughly speaking, if the second variation
of an action functional is positive one knows that the corresponding solution is stable.
More precisely, this means that the critical point is isolated in the space of all solutions
and there does not exist a solution ``nearby''. 
The standard example to illustrate this fact are geodesics in hyperbolic space (which are isolated)
in contrast to geodesics on the sphere (in which case we
can always find a geodesic nearby). This fact also holds true for higher dimensional domains:
Harmonic maps to targets with negative curvature are stable \cite{MR214004}
whereas harmonic maps to spheres
will be unstable in general due to a classic result of Leung \cite{MR673586}.

In quantum field theory several notions of stability exist, however they are usually not connected
to the second variation of the action functional.
Nevertheless, we would like to point out that also for the calculations in physics, in particular in the context
of quantum field theory on an Euclidean space, an understanding of the second
variation should be of great importance.

Assuming that we have a critical point of an action functional in quantum field theory
(which corresponds to a classical solution in physics) one would like to quantize this solution. 
Let us suppose that one can perform this quantization procedure in a rigorous mathematical manner,
then one can think of it as calculating small corrections to the classical solution. 
However, if the classical solution is not stable, any small perturbation of it might change
the classical solution to a nearby solution. 
This might have the effect that one effectively does not calculate the correction
to the classical solution one is looking for but just switches to a classical
solution nearby.

Note that these considerations are only valid if one chooses a Riemannian domain
as action functionals on Lorentzian manifolds are not positive definite in general  
and consequently it does not make much sense to investigate the second variation in this case.

The \(SO(n)\)-nonlinear sigma model is a well-studied 
model in quantum field theory. 
In the physics literature this model is considered both on Riemannian and Lorentzian manifolds
and one often employs the method of Wick rotation to switch between both of them.
It should of course be noted that the method of Wick rotation is not a rigorous
mathematical operation in general.
In mathematical terms the \(SO(n)\)-nonlinear sigma model 
corresponds to harmonic/wave maps with a spherical target of dimension \(n\).

Hence, in the case of a Riemannian domain, our above remarks 
on the instability of harmonic maps to spheres can be applied
to the \(SO(n)\)-nonlinear sigma model.
This model is often used to describe the low-energy behavior of certain particles, 
the so-called \emph{pions}, see \cite[Chapter 19]{MR2148467} and references therein for more 
details on the physics background.
However, one observes that pions are unstable in the sense that they decay into other particles.
This might very well be a consequence of the fact that the second variation of the 
corresponding action functional is not positive.
Consequently, it may not be necessary to perform any kind of quantization procedure
as is done in physics to explain the instability of pions.

Moreover, as we have seen the critical points of the supersymmetric version of the 
nonlinear sigma model will be unstable as well in general 
even if we take into account an additional scalar potential.

\section{Regularity of weak solutions}
In this section we investigate the regularity of weak solution of the system \eqref{euler-lagrange-phi}, \eqref{euler-lagrange-psi}.
Our analysis is very similar to the regularity analysis of the nonlinear sigma model coupled to a gravitino
which was carried out in \cite[Section 2]{MR3798022} 
and to the one for Dirac-harmonic maps with curvature term \cite[Section 4.1]{MR3333092} for \(\dim=2\)
and \cite{MR4018319} for the higher-dimensional case.
The methods used to control the gravitino of \cite{MR3798022} can be modified such that they also apply to our setup.

Let us emphasize again that we are assuming the potential \(V(\phi,\psi)\) to be smooth.
The case of a non-smooth potential can be treated with the help of the methods
developed in \cite{MR3709298} for the nonlinear sigma model coupled to a non-smooth gravitino.

Since we assume to have a compact target manifold the contributions of the potential \(V(\phi,\psi)\)
that do not depend on \(\psi\) are always bounded by a constant.

Motivated by the structure of various potentials that appear in quantum field theory
we will first focus on a potential of the following form
\begin{align}
\label{potential-structure}
V(\phi,\psi)=H(\phi)+G_{\alpha_1\alpha_2\ldots\alpha_s}(\phi)\langle\psi^{\alpha_1},\psi^{\alpha_2}\rangle_{\Sigma M}\times\ldots\times\langle\psi^{\alpha_{s-1}},\psi^{\alpha_{s}}\rangle_{\Sigma M}.
\end{align}
Here, \(H(\phi), G(\phi)\) are smooth functions on \(N\)
and the indices \(0\leq\alpha_i\leq\ldots\leq\alpha_s\) are used to contract the spinors in the potential,
where \(s\) needs to be an even number.

Hence, for the class of potentials from above the following growth assumption holds
\begin{align}
\label{growth-potential}
|V(\phi,\psi)|\leq C_1+C_2|\psi|^s,
\end{align}
for some positive constants \(C_1,C_2\).
In particular, we have
\begin{align}
\label{assumption-growth-potential}
|V_\phi(\phi,\psi)|&\leq C_3+C_4|\psi|^s,\\
\nonumber|V_\psi(\phi,\psi)|&\leq C_5|\psi|^{s-1}
\end{align}
for positive constants \(C_i,i=3,\ldots,5\).

First, we define a weak solution of the system \eqref{euler-lagrange-phi}, \eqref{euler-lagrange-psi}.
\begin{Dfn}
Suppose that the potential \(V(\phi,\psi)\) is of the structure \eqref{potential-structure}.
A weak solution \((\phi,\psi)\) of the system \eqref{euler-lagrange-phi}, \eqref{euler-lagrange-psi}
is a critical point of the action functional \eqref{action-functional}
in the Sobolev space 
\begin{align*}
\chi_t(M,N):=W^{1,2}(M,N)\times W^{1,\frac{4}{3}}(\Sigma M\otimes\phi^\ast TN)\times L^t(\Sigma M\otimes\phi^\ast TN),
\end{align*}
where \(t=4\) if \(0<s\leq 4\) and \(t=s\) if \(s>4\). 
\end{Dfn}

\begin{Bem}
The way we have defined the space \(\chi_t(M,N)\) is inspired from the analysis of Dirac-harmonic maps
which corresponds to solutions of \eqref{euler-lagrange-phi}, \eqref{euler-lagrange-psi} with \(V(\phi,\psi)=0\).
This fits naturally to how we have defined weak solutions above
as the action functional \eqref{action-functional} is finite whenever \((\phi,\psi)\in\chi_t(M,N)\).
For a potential with stronger growth as in \eqref{growth-potential} we 
will have to adjust the space \(\chi_t(M,N)\) in such a way that \eqref{action-functional} is finite.
\end{Bem}

In order to state the main regularity result we recall the following:
\begin{Dfn}
For a given open subset \(D\subset\R^m\) and \(1\leq p<\infty\), \(0<\lambda\leq m\)
the Morrey space \(M^{p,\lambda}(D)\) is defined as follows
\begin{align*}
M^{p,\lambda}(D):=\left\{f\in L^p(D)\mid \|f\|^p_{M^{p,\lambda}(D)}:=\sup_{B_r\subset D}\{r^{\lambda-m}\int_{B_r}|f|^p(y)dy\}<\infty\right\}.
\end{align*}
\end{Dfn}
Note that \(M^{p,m}(D)=L^p(D)\).

The following smallness condition will be important in the 
formulation of our main regularity result:
\begin{align}
\label{assumption-regularity-theorem}
\|d\phi\|_{M^{2,2}(U)}+\|\psi\|_{M^{4,2}(U)}\leq\epsilon
\end{align}
Here, \(\epsilon>0\) is a small number depending on \((M,g), (N,h)\) 
and \(U\) an open subset of \(M\).

We will prove the following regularity result for weak solutions of the Dirac-harmonic map
with potential system:
\begin{Satz}
\label{theorem-regularity}
Let \(m\geq 2\) and suppose that the potential \(V(\phi,\psi)\) is smooth and has the structure \eqref{potential-structure}.
We need to make the following case distinction:
\begin{enumerate}
 \item If \(1<s\leq 4\) suppose that \eqref{assumption-regularity-theorem} holds.
 \item If \(s>4\) suppose that \eqref{assumption-regularity-theorem} holds and also \(\psi\in{M^{2,2s-2}(U)}\).
\end{enumerate}
Then a weak solution \((\phi,\psi)\in\chi_t(M,N)\) is smooth in \(U\), where \(U\) is an open subset of \(M\).
\end{Satz}

\begin{Bem}
\begin{enumerate}
\item 
 For \(\dim M=2\) the smallness condition \eqref{assumption-regularity-theorem} can easily be achieved.
In this case it reads \(\|d\phi\|_{L^2(U)}+\|\psi\|_{L^4(U)}\leq\epsilon\)
and can be satisfied by choosing \(U\) sufficiently small.
\item Due to the Sobolev embedding theorem we have \(\|\psi\|_{L^\frac{m}{\frac{3}{4}m-1}(U)}\leq C\|\nabla\psi\|_{L^\frac{4}{3}(U)}\)
such that for \(\dim M=2\) we do not have to require that \(\psi\in L^4(\Sigma M\otimes\phi^\ast TN)\) in the definition of the space \(\chi_t(M,N)\)
if \(0<s\leq 4\).
The Sobolev embedding theorem does not help in the case of \(m\geq 3\) as we have \(\frac{m}{\frac{3}{4}m-1}<4\).
\item Note that Theorem \ref{theorem-regularity} and its proof can easily be modified such that they also hold for non-integer values of \(s\).
\end{enumerate}
\end{Bem}

In order to prove Theorem \ref{theorem-regularity} it is favorable to apply the embedding theorem of Nash and to isometrically
embed the target manifold \(N\) into some \(\R^q\) of sufficiently large dimension \(q\).
We denote this isometric embedding by \(\iota\colon N\to\R^q\) which we assume to be smooth.
Moreover, let \(\phi':=\iota\circ\phi\colon M\to\R^q\) and \(\psi'=d\iota(\psi)\) be 
the corresponding pushforward of the vector spinor. As \(\psi\in\Gamma(\Sigma M\otimes\phi^\ast TN)\)
behaves like a tangent vector on \(N\) we have to use the differential of \(\iota\)
to push it to the ambient space \(\R^q\), hence we have \(\psi'=d\iota(\psi)\in\Gamma(\Sigma M\otimes\R^q)\).

Set \(n:=\dim N\) and 
let \(u^\alpha,1\leq\alpha\leq q\) be global coordinates on the ambient space \(\R^q\).
Moreover, let \(\nu_\theta,\theta=n+1,\ldots,q\) be an orthonormal frame of the submanifold \(\iota(\phi)\).
In this setup \(\phi'\) becomes a vector-valued function and \(\psi'\)
can be thought of as a vector of standard spinors \(\psi'^\alpha\in\Gamma(\Sigma M),1\leq\alpha\leq q\)
that are constrained as follows
\begin{align*}
\sum_{\alpha}\nu^\alpha_\theta\psi'^\alpha=0,\qquad n+1\leq\theta\leq q.
\end{align*}

This allows us to give the extrinsic version of the 
Euler-Lagrange equations \eqref{euler-lagrange-phi}, \eqref{euler-lagrange-psi}. 

In the following we will omit the \('\) and again write \(\phi,\psi\) for the extrinsic 
version of the Euler-Lagrange equations to shorten the notation.

\begin{Lem}
The extrinsic version of the Euler-Lagrange equations \eqref{euler-lagrange-phi}, \eqref{euler-lagrange-psi} 
for \(\phi\colon M\to\R^q\) and \(\psi\in\Gamma(\Sigma M\otimes\R^q)\) is given by the system
\begin{align}
\label{extrinsic-phi}
\Delta\phi^\alpha&=(\omega_{i}^{\alpha\beta}+F_i^{\alpha\beta})\frac{\partial\phi^\beta}{\partial x^i}-(V_\phi(\phi,\psi))^\alpha, \\
\label{extrinsic-psi}
\p\psi^\alpha&=-\nabla\phi^\delta\cdot\psi^\beta\frac{\partial\nu_l^\beta}{\partial u^\delta}\nu_l^\alpha+(V_\psi(\phi,\psi))^\alpha,
\end{align}
where \(1\leq\alpha\leq q\) and
\begin{align*}
\omega_{i}^{\alpha\beta}&=
\frac{\partial\phi^\gamma}{\partial x^i}\frac{\partial\nu_l^\beta}{\partial u^\gamma}\nu_l^\alpha-\frac{\partial\phi^\gamma}{\partial x^i}\frac{\partial\nu_l^\alpha}{\partial u^\gamma}\nu_l^\beta
=-\omega_{i}^{\alpha\beta},\\
F_i^{\alpha\beta}&=\langle\psi^\gamma,e_i\cdot\psi^\delta\rangle
\bigg(\big(\frac{\partial\nu_l}{\partial u^\delta}\big)^{\top,\beta}\big(\frac{\partial\nu_l}{\partial u^\gamma}\big)^{\top,\alpha}
-\big(\frac{\partial\nu_l}{\partial u^\delta}\big)^{\top,\alpha}\big(\frac{\partial\nu_l}{\partial u^\gamma}\big)^{\top,\beta}
\bigg)
=-F_i^{\beta\alpha}.
\end{align*}
Here, \(\top\) denotes the projection map \(\top\colon\R^q\to TN\).
\end{Lem}
\begin{proof}
We will only treat the terms involving the potential \(V(\phi,\psi)\),
for the remaining terms see the detailed discussion in \cite[Section 4.3]{MR3772035}.

Suppose that \(\tilde V_\phi(\phi',\psi')\) is an extension of \(V_\phi(\phi,\psi)\) to the ambient space \(\R^q\).
Then, using the projection \(\top\) to the tangent space \(T_yN\) at the point \(y\in N\) 
we can set 
\begin{align*}
(\tilde V_\phi(\phi',\psi'))^\top=V_\phi(\phi,\psi).
\end{align*}
We can also extend \(V_\psi(\phi,\psi)\) to the ambient space, denote this extension by \(\tilde V_\psi(\phi,\psi)\).
By the same argument as before we set
\begin{align*}
(\tilde V_\psi(\phi',\psi'))^\top=V_\psi(\phi,\psi)
\end{align*}
completing the proof.
\end{proof}

Abbreviating
\begin{align*}
\Omega^{\alpha\beta}_i=\omega_{i}^{\alpha\beta}+F_i^{\alpha\beta},\qquad A_{\alpha\beta}=-\nabla\phi^\delta\frac{\partial\nu_l^\beta}{\partial u^\delta}\nu_l^\alpha
\end{align*}
the extrinsic version of the Euler-Lagrange equations \eqref{extrinsic-phi}, \eqref{extrinsic-psi} can be written in the compact form
\begin{align}
\label{euler-lagrange-antisymmetric}
\Delta\phi&=\Omega\cdot\nabla\phi-V_\phi(\phi,\psi),\\
\nonumber\p\psi&=A\cdot\psi+V_\psi(\phi,\psi)
\end{align}
with the important feature that \(\Omega\) is antisymmetric, that is \(\Omega^{\alpha\beta}=-\Omega^{\beta\alpha}\).
Note that we have \(|\Omega|^2\leq C(|d\phi|^2+|\psi|^4)\) and \(|A|\leq C|d\phi|\).

At this point we are ready to apply tools from elliptic regularity that will
help to improve the integrability of the system \eqref{euler-lagrange-antisymmetric}.

On the one hand we will make use of
the following Lemma, which was proven in \cite[Lemma 6.1]{MR3772035}:

\begin{Lem}
\label{lemma-spinor}
Let \(B_1\subset\R^m\), \(m\geq 2\) and let \(4<p<\infty\). Consider a weak solution \(\varphi\in M^{4,2}(B_1,\R^L\otimes\R^q)\) of
\begin{align}
\p\varphi^\alpha=A_{\alpha\beta}\cdot\varphi^\beta+B^\alpha,\qquad 1\leq\alpha\leq q,
\end{align}
where \(A\in M^{2,2}(B_1,\mathfrak{gl}(L,\R)\otimes\mathfrak{gl}(q,\R))\) and \(B\in M^{2,2}(B_1,\R^L\otimes\R^q)\).
Then, there exists \(\epsilon_0=\epsilon_0(m,p)>0\) such that if
\begin{align*}
\|A\|_{M^{2,2}(B_1)}\leq\epsilon_0,
\end{align*}
we have \(\varphi\in L^p_{loc}(B_1)\).
Moreover, for any \(U\subset B_1\) there holds the estimate
\begin{align}
\|\varphi\|_{L^p(U)}\leq C(m,p,U)(\|\varphi\|_{M^{4,2}(B_1)}+\|B\|_{M^{2,2}(B_1)})
\end{align}
for some \(C(m,p,U)>0\).
\end{Lem}

On the other hand, in order to improve the integrability of the map \(\phi\), we employ the following regularity result,
which was obtained in \cite[Theorem 1.2]{MR3251921}: 
\begin{Satz}
\label{theorem-sharp-hd}
Let \(B_1\subset\R^m\), \(m\geq 2\) and let \(u\in W^{1,2}(B_1,\R^q)\) with \(\nabla u\in M^{2,2}(B_1,\R^m\otimes\R^q)\),
\(\Omega\in M^{2,2}(B_1,\mathfrak{so}(q)\otimes\R^m)\) and \(f\in L^p(B_1)\) for \(\frac{m}{2}<p<m\),
be a weak solution of 
\begin{align*}
-\Delta u=\Omega\cdot\nabla u+f.
\end{align*}
Then, for any \(U\subset B_1\) there exists \(\epsilon_1=\epsilon_1(m,p,q,U)\) such that whenever \(\|\Omega\|_{M^{2,2}}\leq\epsilon_1\)
the following estimate holds
\begin{align}
\|\nabla^2u\|_{M^{\frac{2p}{m},2}(U)}+\|\nabla u\|_{M^{\frac{2p}{m-p},2}(U)}\leq C(\|u\|_{L^1(B_1)}+\|f\|_{L^p(B_1)}).
\end{align}
\end{Satz}
Note that in \cite{MR3251921} a different convention in the definition of Morrey spaces is used
as in this article.

\begin{Bem}
We can think of Lemma \ref{lemma-spinor} as providing a regularity result for 
Dirac equations coupled to a potential, whereas 
Theorem \ref{theorem-sharp-hd} constitutes a regularity result for Laplace-tpye
equations allowing a gradient term on the right hand side with an antisymmetric structure
and also a potential term. While both results make use of the ellipticity of
the Dirac and the Laplace operator it becomes clear that the smoothing effects
of the Laplacian are stronger: In Lemma \ref{lemma-spinor} we need to demand 
that the potential lies in some Morrey space whereas in Theorem \ref{theorem-sharp-hd}
it is enough that the potential is integrable in \(L^p\) for specific values of \(p\).
\end{Bem}

We are now ready to give the proof of Theorem \ref{theorem-regularity}.

\begin{proof}[Proof of Theorem \ref{theorem-regularity}]
First, we will improve the regularity of the spinor by making use of Lemma \ref{lemma-spinor}.
Thus, we use a local trivialization of \(\Sigma M\) over \(B_1\) and make the following case distinction:
\begin{enumerate}
 \item If \(1<s<4\) the only case to be considered is \(s=2\).
  Then, we apply Lemma \ref{lemma-spinor} with \(|B|\leq C|\psi|\).
  Note that \(\|\psi\|_{M^{2,2}(U)}\leq C\|\psi\|_{M^{4,2}(U)}\leq C\epsilon\)
  due to the Hölder inequality in Morrey spaces.
 \item If \(s=4\) we can write the equation for the spinor \eqref{extrinsic-psi} in the form
  \begin{align*}
   \p\varphi^\alpha=A_{\alpha\beta}\cdot\varphi^\beta,\qquad 1\leq\alpha\leq q,
  \end{align*}
   where \(|A|\leq C(|d\phi|+|\psi|^2)\) such that we can apply Lemma \ref{lemma-spinor} with \(B=0\).
 \item If \(s>4\) we again apply Lemma \ref{lemma-spinor} with \(|B|\leq C|\psi|^{s-1}\).
\end{enumerate}
Due to the assumptions made in Theorem \ref{theorem-regularity} we get from Lemma \ref{lemma-spinor}
that \(\psi\in L^p_{loc}\) for any \(p\in[1,\infty)\).

At this point we are ready to improve the integrability of the map \(\phi\).
Using the growth condition for the potential \eqref{growth-potential}
we can apply Theorem \ref{theorem-sharp-hd}. Setting \(p=m-\delta\) we calculate
\begin{align*}
\frac{2p}{m-p}=\frac{2(m-\delta)}{m-(m-\delta)}=\frac{2m}{\delta}-2
\end{align*}
and choosing \(\delta>0\) sufficiently small 
we obtain the continuity of \(\phi\) from the Morrey Lemma.
From this point on Theorem \ref{theorem-regularity} follows by a standard bootstrap procedure.
\end{proof}

Let us shortly discuss how Theorem \ref{theorem-regularity} needs to be modified if we 
would consider an arbitrary potential \(V(\phi,\psi)\) that does not have the structure \eqref{potential-structure}.
In the latter case we would need to consider the function space
\begin{align*}
\tilde\chi(M,N):=W^{1,2}(M,N)\times W^{1,\frac{4}{3}}(\Sigma M\otimes\phi^\ast TN)\times L^4(\Sigma M\otimes\phi^\ast TN)\times T(\Sigma M\otimes\phi^\ast TN),
\end{align*}
where the space \(T(\Sigma M\otimes\phi^\ast TN)\) needs to be chosen such that \(S_P(\phi,\psi)\) is finite for the potential under consideration.

\begin{Satz}
Let \(m\geq 2\), \(U\subset M\) open and suppose that the potential \(V(\phi,\psi)\) is smooth and satisfies \(\|V_\psi(\phi,\psi)\|_{M^{2,2}}(U)<\infty\).
Then, there exists \(\epsilon>0\) such that if
\begin{align*}
\|d\phi\|_{M^{2,2}(U)}+\|\psi\|_{M^{4,2}(U)}\leq\epsilon,
\end{align*}
where \(\epsilon>0\) is a small number depending on \((M,g), (N,h)\), 
a weak solution \((\phi,\psi)\in\tilde\chi(M,N)\) is smooth in \(U\).
\end{Satz}
\begin{proof}
This again follows from Lemma \ref{lemma-spinor} and Theorem \ref{theorem-sharp-hd}.
\end{proof}

\section{Explicit potentials from quantum field theory}
In this section we will apply our previous analysis for an arbitrary potential
to some specific potentials which are studied in quantum field theory.
We will shed light on the particular mathematical advantages and drawbacks
of the potentials used in theoretical physics.

\subsection{Dirac-harmonic maps with curvature term}
Dirac-harmonic maps with curvature term arise as critical points of 
\eqref{action-functional} by choosing the following potential
\begin{align*}
V_1(\phi,\psi):=\frac{1}{12}\langle R^N(\phi)(\psi,\psi)\psi,\psi\rangle.
\end{align*}
Here, the spinors are contracted in the following way
\[
\langle R^N(\psi,\psi)\psi,\psi\rangle=R_{\alpha\beta\gamma\delta}\langle\psi^\alpha,\psi^\gamma\rangle\langle\psi^\beta,\psi^\delta\rangle,
\]
which ensures that the action functional is real-valued.
This particular potential was first studied in the mathematics literature in \cite{MR2370260}.

In the physics literature the prefactor in front of the potential is chosen
such that it ensures the invariance of the action functional under supersymmetry transformations \cite{MR626710}.
However, the mathematical analysis initiated in \cite{MR2370260} applies to any non-zero constant
in front of the curvature term.

The equations for Dirac-harmonic maps with curvature term are the following:
\begin{align}
\label{phi-dhc}\tau(\phi)=&\frac{1}{2}R^N(\psi,e_i\cdot\psi)d\phi(e_i)
-\frac{1}{12}\langle(\nabla R^N)^\sharp(\psi,\psi)\psi,\psi\rangle, \\
\label{dhc-psi}\D\psi=&\frac{1}{3}R^N(\psi,\psi)\psi,
\end{align}
where \(\sharp\colon\phi^\ast T^\ast N\to\phi^\ast TN\) represents the musical isomorphism.

For the domain being a closed surface it was shown in \cite{MR3333092} that a weak solution \((\phi,\psi)\in\chi_{t}(M,N)\) of \eqref{euler-lagrange-phi}, \eqref{euler-lagrange-psi} with \(t=4\)
is smooth which was later extended to higher-dimensional domains in \cite{MR4018319}.

If the domain manifold is two-dimensional the action functional for 
Dirac-harmonic maps with curvature term is invariant under conformal transformations
such that in this case the solutions of the resulting Euler-Lagrange equations share nice properties,
one of them being the possibility of removing isolated point singularities 
if a certain energy is finite \cite{MR3558358}.

Note that the smoothness of weak solutions of \eqref{phi-dhc}, \eqref{dhc-psi}
can also also be established with the help of Theorem \ref{theorem-regularity},
where we would consider the case \(s=4\).

However, it is not hard to realize that the action functional for Dirac-harmonic maps
with curvature term is still unbounded from below meaning that the curvature term
does not help in curing the analytic problems of the Dirac-harmonic map action functional.

\subsection{A potential from supersymmetric quantum field theory}
As already explained in the introduction the action functionals of quantum field theory 
are built by the requirement that they are invariant under certain symmetry operations.
It was shown in \cite{MR719813}, see also \cite[Theorem 3.82]{MR1701600},
that the following potential 
\begin{align*}
V_2(\phi,\psi)=\frac{1}{2}|\nabla W(\phi)|^2+\frac{1}{2}\hess W(\phi)(\psi,\psi)-\frac{1}{12}\langle R^N(\phi)(\psi,\psi)\psi,\psi\rangle
\end{align*}
preserves the invariance under supersymmetry transformations
making it a distinct potential for the supersymmetric nonlinear sigma model.
Here, \(W\colon N\to\R\) is a smooth function.
On the other hand it is not hard to realize that this particular potential
does not preserve the invariance under conformal transformations
on a two-dimensional domain manifold.

The resulting Euler-Lagrange equations acquire the following form
\begin{align}
\label{phi-pot}\tau(\phi)=&\frac{1}{2}R^N(\psi,e_i\cdot\psi)d\phi(e_i)
+\frac{1}{12}\langle(\nabla R^N)^\sharp(\psi,\psi)\psi,\psi\rangle
-\frac{1}{2}\nabla|\nabla W(\phi)|^2-\frac{1}{2}\nabla\hess W(\phi)(\psi,\psi), \\
\label{dhc-pot}\D\psi=&-\frac{1}{3}R^N(\psi,\psi)\psi+\nabla_{(\cdot)}dW(\psi).
\end{align}

In the case of a compact target \(N\) the potential \(V_2(\phi,\psi)\) satisfies
the following growth conditions
\begin{align*}
|V_2(\phi,\psi)|&\leq C_1+C_2|\psi|^2+C_3|\psi|^4,\\
|V_{2\phi}(\phi,\psi)|&\leq C_4+C_5|\psi|^2+C_6|\psi|^4,\\
|V_{2\psi}(\phi,\psi)|&\leq C_7|\psi|+C_8|\psi|^3
\end{align*}
for positive constants \(C_i,i=1,\ldots 8\).
In order to apply our main regularity result, that is Theorem \ref{theorem-regularity},
we need to satisfy the assumptions of case (1) with \(s=2\) and \(s=4\).
Hence, a weak solution \((\phi,\psi)\in\chi_t(M,N)\) with \(t=4\) is smooth
if the smallness condition \eqref{assumption-regularity-theorem} 
holds.

On the other hand, it can again easily be checked that the potential \(V_2(\phi,\psi)\)
will not be helpful in achieving the positivity of the action functional for Dirac-harmonic maps
or the positivity of the second variation although it respects the invariance under 
supersymmetry transformations in quantum field theory.

\subsection{Some further potentials}
Let us give a short overview on some further potentials that fit into 
the framework of this article.

\begin{enumerate}
 \item In quantum field theory one often includes a term in the action functional that models the mass of a fermion.
 From a mathematical point of view this is reflected in considering a potential of the form 
 \begin{align*}
 V_3(\phi,\psi)=\frac{\lambda}{2}|\psi|^2, \qquad \lambda>0.
 \end{align*}
However, this particular potential is again not helpful to achieve a lower bound on the action functional \eqref{action-functional}.
In addition, including the mass term destroys the conformal invariance of the action functional.
The regularity of Dirac-harmonic maps coupled to a similar potential was studied in \cite{MR3018163},
where the Ricci curvature of the target manifold is used to define a potential that is quadratic in the spinors.		

The resulting Euler-Lagrange equations are 
\begin{align}
\label{phi-mass}\tau(\phi)=&\frac{1}{2}R^N(\psi,e_i\cdot\psi)d\phi(e_i), \\
\label{dhc-mass}\D\psi=&\lambda\psi.
\end{align}
Note that the second equation states that we should find an eigenspinor of the twisted Dirac operator \(\D\)
with positive eigenvalue \(\lambda\).

The potential \(V_3(\phi,\psi)\) satisfies
the following growth conditions
\begin{align*}
|V_3(\phi,\psi)|&\leq C_1|\psi|^2,\\
|V_{3\psi}(\phi,\psi)|&\leq C_2|\psi|
\end{align*}
for positive constants \(C_i,i=1,2\).
To apply Theorem \ref{theorem-regularity} again
we need to satisfy the assumptions of case (1) with \(s=2\).
Thus, as for the previous potential, a weak solution \((\phi,\psi)\in\chi_t(M,N)\) with \(t=4\) is smooth
if the smallness condition \eqref{assumption-regularity-theorem} 
holds.

On the other hand, we realize again that this specific potential cannot remove the indefinite
character of the action functional for Dirac-harmonic maps.

\item From a purely mathematical perspective a potential of the form
\begin{align*}
V_4(\phi,\psi)=V(\phi)-\exp(|\psi|^2)
\end{align*}
could be beneficial. This particular potential
will not respect the conformal invariance on a two-dimensional domain
and also the regularity theory established in Theorem \ref{theorem-regularity}
does not apply in this case.
On the other hand, the potential \(V_4(\phi,\psi)\) may be helpful in order to 
establish the positivity of the action functional \eqref{action-functional}
and also its second variation \eqref{second-variation}.

Inserting \(V_4(\phi,\psi)\) into the condition for the positivity of the action functional \eqref{positivity-action}
we obtain 
\begin{align*}
\sum_{J=-\infty}^\infty|\alpha_J|^2\lambda_J-2\int_MV(\phi)\dv+2\int_M\exp(|\psi|^2)\dv\geq 0.
\end{align*}
Although the potential \(V_4(\phi,\psi)\) contains an exponential factor, we still cannot
ensure that it can compensate the negative eigenvalues of the Dirac-term.

Moreover, inserting \(V_4(\phi,\psi)\) into the conditions for 
the second variation to be positive, which are \eqref{second-variation-positivity-a} and \eqref{second-variation-positivity-b},
we obtain the requirements
\begin{align*}
-\iota(\eta,\eta)\big(V_{\phi\phi}(\phi)\big)\geq C|\eta|^2e(\phi,\psi)
\end{align*}
and
\begin{align*}
\int_M\big(&\langle\xi,\D\xi\rangle-c_6|\xi|^2|\psi|^2
+2|\xi|^2\exp(|\psi|^2)+4\exp(|\psi|^2)|\langle\psi,\xi\rangle|^2\big)\dv\geq 0.
\end{align*}
While the first condition can be satisfied on a non-compact target manifold \(N\)
if \(e(\phi,\psi)\) is finite we still may encounter problems in controlling the Dirac-term
in the second inequality.

\item In the case of a vanishing spinor, that is \(\psi=0\), the action functional \eqref{action-functional}
reduces to the one for harmonic maps with potential introduced in \cite{MR1433176}.

The corresponding Euler-Lagrange equation is
\begin{align}
\label{harmonic-pot}\tau(\phi)=-V_\phi(\phi).
\end{align}

For harmonic maps with potential we obtain the following variant of Theorem \ref{theorem-regularity}
generalizing \cite[Theorem 1.1]{MR3884770} to higher dimensions:

\begin{Satz}
\label{theorem-regularity-harmonic-potential}
Let \(m\geq 2\) and suppose that the potential \(V(\phi)\) is smooth.
Then there exists \(\epsilon>0\) depending on \((M,g), (N,h)\)
such that if \(\phi\in W^{1,2}(M,N)\)
is a weak solution of \eqref{harmonic-pot} satisfying
\begin{align}
\label{assumption-regularity-theorem-harmonic-potential}
\|d\phi\|_{M^{2,2}(U)}\leq\epsilon
\end{align}
for some open subset \(U\subset M\) then \(\phi\) is smooth in \(U\).
\end{Satz}
\begin{proof}
This follows directly from Theorem \ref{theorem-sharp-hd} and a standard bootstrap.
\end{proof}

\end{enumerate}

\section{Declarations}
\textbf{Conflict of interest statement}:
The corresponding author states that there is no conflict of interest. \\

\textbf{Code availability statement}: Not applicable. \\

\textbf{Data Availibility Statement}: Data sharing not applicable to this article as no datasets were
generated or analysed during the current study.

\bibliographystyle{plain}
\bibliography{mybib}

\def\cprime{$'$}
\begin{thebibliography}{10}

\bibitem{MR626710}
Luis Alvarez-Gaum\'{e} and Daniel~Z. Freedman.
\newblock Geometrical structure and ultraviolet finiteness in the
  supersymmetric {$\sigma $}-model.
\newblock {\em Comm. Math. Phys.}, 80(3):443--451, 1981.

\bibitem{MR719813}
Luis Alvarez-Gaum\'{e} and Daniel~Z. Freedman.
\newblock Potentials for the supersymmetric nonlinear {$\sigma $}-model.
\newblock {\em Comm. Math. Phys.}, 91(1):87--101, 1983.

\bibitem{MR3070562}
Bernd Ammann and Nicolas Ginoux.
\newblock Dirac-harmonic maps from index theory.
\newblock {\em Calc. Var. Partial Differential Equations}, 47(3-4):739--762,
  2013.

\bibitem{MR3305429}
Volker Branding.
\newblock Magnetic {D}irac-harmonic maps.
\newblock {\em Anal. Math. Phys.}, 5(1):23--37, 2015.

\bibitem{MR3333092}
Volker Branding.
\newblock Some aspects of {D}irac-harmonic maps with curvature term.
\newblock {\em Differential Geom. Appl.}, 40:1--13, 2015.

\bibitem{MR3493217}
Volker Branding.
\newblock Dirac-harmonic maps with torsion.
\newblock {\em Commun. Contemp. Math.}, 18(4):1550064, 19, 2016.

\bibitem{MR3558358}
Volker Branding.
\newblock Energy estimates for the supersymmetric nonlinear sigma model and
  applications.
\newblock {\em Potential Anal.}, 45(4):737--754, 2016.

\bibitem{MR3735550}
Volker Branding.
\newblock On conservation laws for the supersymmetric sigma model.
\newblock {\em Results Math.}, 72(4):2181--2201, 2017.

\bibitem{MR3673634}
Volker Branding.
\newblock Some remarks on energy inequalities for harmonic maps with potential.
\newblock {\em Arch. Math. (Basel)}, 109(2):151--165, 2017.

\bibitem{MR3884770}
Volker Branding.
\newblock A global weak solution to the full bosonic string heat flow.
\newblock {\em J. Evol. Equ.}, 18(4):1819--1841, 2018.

\bibitem{MR3830780}
Volker Branding.
\newblock A note on twisted {D}irac operators on closed surfaces.
\newblock {\em Differential Geom. Appl.}, 60:54--65, 2018.

\bibitem{MR3886921}
Volker Branding.
\newblock A vanishing result for the supersymmetric nonlinear sigma model in
  higher dimensions.
\newblock {\em J. Geom. Phys.}, 134:1--10, 2018.

\bibitem{MR3917346}
Volker Branding.
\newblock Energy methods for {D}irac-type equations in two-dimensional
  {M}inkowski space.
\newblock {\em Lett. Math. Phys.}, 109(2):295--325, 2019.

\bibitem{MR4034775}
Volker Branding.
\newblock Nonlinear {D}irac {E}quations, {M}onotonicity {F}ormulas and
  {L}iouville {T}heorems.
\newblock {\em Comm. Math. Phys.}, 372(3):733--767, 2019.

\bibitem{MR3830277}
Volker Branding and Klaus Kr\"{o}ncke.
\newblock Global existence of {D}irac-wave maps with curvature term on
  expanding spacetimes.
\newblock {\em Calc. Var. Partial Differential Equations}, 57(5):Art. 119, 30,
  2018.

\bibitem{MR2370260}
Q.~Chen, J.~Jost, and G.~Wang.
\newblock Liouville theorems for {D}irac-harmonic maps.
\newblock {\em J. Math. Phys.}, 48(11):113517, 13, 2007.

\bibitem{MR2262709}
Qun Chen, J{\"u}rgen Jost, Jiayu Li, and Guofang Wang.
\newblock Dirac-harmonic maps.
\newblock {\em Math. Z.}, 254(2):409--432, 2006.

\bibitem{MR1701600}
Pierre Deligne and Daniel~S. Freed.
\newblock Supersolutions.
\newblock In {\em Quantum fields and strings: a course for mathematicians,
  {V}ol. 1, 2 ({P}rinceton, {NJ}, 1996/1997)}, pages 227--355. Amer. Math.
  Soc., Providence, RI, 1999.

\bibitem{MR1433176}
Ali Fardoun and Andrea Ratto.
\newblock Harmonic maps with potential.
\newblock {\em Calc. Var. Partial Differential Equations}, 5(2):183--197, 1997.

\bibitem{MR1800592}
Ali Fardoun, Andrea Ratto, and Rachid Regbaoui.
\newblock On the heat flow for harmonic maps with potential.
\newblock {\em Ann. Global Anal. Geom.}, 18(6):555--567, 2000.

\bibitem{MR2138082}
Xiaoli Han.
\newblock Dirac-wave maps.
\newblock {\em Calc. Var. Partial Differential Equations}, 23(2):193--204,
  2005.

\bibitem{MR214004}
Philip Hartman.
\newblock On homotopic harmonic maps.
\newblock {\em Canadian J. Math.}, 19:673--687, 1967.

\bibitem{MR2860404}
Takeshi Isobe.
\newblock On the existence of nonlinear {D}irac-geodesics on compact manifolds.
\newblock {\em Calc. Var. Partial Differential Equations}, 43(1-2):83--121,
  2012.

\bibitem{MR3978518}
Takeshi Isobe.
\newblock On the multiple existence of superquadratic {D}irac-harmonic maps
  into flat tori.
\newblock {\em Calc. Var. Partial Differential Equations}, 58(4):Art. 126, 41,
  2019.

\bibitem{isobe2017morse}
Takeshi Isobe and Ali Maalaoui.
\newblock Morse-{F}loer theory for super-quadratic {D}irac-geodesics.
\newblock {\em arXiv preprint arXiv:1712.08960}, 2017.

\bibitem{MR3772035}
J\"{u}rgen Jost, Enno Ke{\ss}ler, J\"{u}rgen Tolksdorf, Ruijun Wu, and Miaomiao
  Zhu.
\newblock Regularity of solutions of the nonlinear sigma model with gravitino.
\newblock {\em Comm. Math. Phys.}, 358(1):171--197, 2018.

\bibitem{MR3913850}
J\"{u}rgen Jost, Enno Ke{\ss}ler, J\"{u}rgen Tolksdorf, Ruijun Wu, and Miaomiao
  Zhu.
\newblock From harmonic maps to the nonlinear supersymmetric sigma model of
  quantum field theory: at the interface of theoretical physics, {R}iemannian
  geometry, and nonlinear analysis.
\newblock {\em Vietnam J. Math.}, 47(1):39--67, 2019.

\bibitem{MR3724759}
J\"{u}rgen Jost, Lei Liu, and Miaomiao Zhu.
\newblock A global weak solution of the {D}irac-harmonic map flow.
\newblock {\em Ann. Inst. H. Poincar\'{e} Anal. Non Lin\'{e}aire},
  34(7):1851--1882, 2017.

\bibitem{MR4018319}
J\"{u}rgen Jost, Lei Liu, and Miaomiao Zhu.
\newblock Regularity of {D}irac-harmonic maps with {$\lambda$}-curvature term
  in higher dimensions.
\newblock {\em Calc. Var. Partial Differential Equations}, 58(6):Paper No. 187,
  24, 2019.

\bibitem{MR3709298}
J\"{u}rgen Jost, Ruijun Wu, and Miaomiao Zhu.
\newblock Coarse regularity of solutions to a nonlinear sigma-model with
  {$L^p$} gravitino.
\newblock {\em Calc. Var. Partial Differential Equations}, 56(6):Art. 154, 17,
  2017.

\bibitem{MR3798022}
J\"{u}rgen Jost, Ruijun Wu, and Miaomiao Zhu.
\newblock Partial regularity for a nonlinear sigma model with gravitino in
  higher dimensions.
\newblock {\em Calc. Var. Partial Differential Equations}, 57(3):Art. 85, 17,
  2018.

\bibitem{MR4259181}
J\"{u}rgen Jost and Jingyong Zhu.
\newblock {$\alpha $}-{D}irac-harmonic maps from closed surfaces.
\newblock {\em Calc. Var. Partial Differential Equations}, 60(3):Paper No. 111,
  41, 2021.

\bibitem{MR4310166}
J\"{u}rgen Jost and Jingyong Zhu.
\newblock Existence of ({D}irac-)harmonic maps from degenerating (spin)
  surfaces.
\newblock {\em J. Geom. Anal.}, 31(11):11165--11189, 2021.

\bibitem{MR4232501}
J\"{u}rgen Jost and Jingyong Zhu.
\newblock Short-time existence of the {$\alpha$}-{D}irac-harmonic map flow and
  applications.
\newblock {\em Comm. Partial Differential Equations}, 46(3):442--469, 2021.

\bibitem{MR1031992}
H.~Blaine Lawson, Jr. and Marie-Louise Michelsohn.
\newblock {\em Spin geometry}, volume~38 of {\em Princeton Mathematical
  Series}.
\newblock Princeton University Press, Princeton, NJ, 1989.

\bibitem{MR673586}
Pui~Fai Leung.
\newblock On the stability of harmonic maps.
\newblock In {\em Harmonic maps ({N}ew {O}rleans, {L}a., 1980)}, volume 949 of
  {\em Lecture Notes in Math.}, pages 122--129. Springer, Berlin-New York,
  1982.

\bibitem{MR604040}
J.~Sacks and K.~Uhlenbeck.
\newblock The existence of minimal immersions of {$2$}-spheres.
\newblock {\em Ann. of Math. (2)}, 113(1):1--24, 1981.

\bibitem{MR3251921}
Ben Sharp.
\newblock Higher integrability for solutions to a system of critical elliptic
  {PDE}.
\newblock {\em Methods Appl. Anal.}, 21(2):221--240, 2014.

\bibitem{MR2544729}
Changyou Wang and Deliang Xu.
\newblock Regularity of {D}irac-harmonic maps.
\newblock {\em Int. Math. Res. Not. IMRN}, (20):3759--3792, 2009.

\bibitem{MR2148467}
Steven Weinberg.
\newblock {\em The quantum theory of fields. {V}ol. {II}}.
\newblock Cambridge University Press, Cambridge, 2005.
\newblock Modern applications.

\bibitem{MR3018163}
Deliang Xu and Zhengxiang Chen.
\newblock Regularity for {D}irac-harmonic map with {R}icci type spinor
  potential.
\newblock {\em Calc. Var. Partial Differential Equations}, 46(3-4):571--590,
  2013.

\end{thebibliography}
\end{document}